\input amstex\documentstyle{amsppt}  
\pagewidth{12.5cm}\pageheight{19cm}\magnification\magstep1  
\topmatter
\title Two partitions of a flag manifold\endtitle
\author G. Lusztig\endauthor
\address{Department of Mathematics, M.I.T., Cambridge, MA 02139}\endaddress
\thanks{Supported by NSF grant DMS-1855773}\endthanks
\endtopmatter   
\document

\define\pos{\text{\rm pos}}
\define\Irr{\text{\rm Irr}}

\define\si{\sim}

\define\sqc{\sqcup}

\define\op{\oplus}
   
\define\part{\partial}
\define\emp{\emptyset}

\define\n{\notin}

\define\m{\mapsto}
\define\do{\dots}

\define\bsl{\backslash}

\define\lra{\leftrightarrow}

\define\sub{\subset}    

\define\T{\times}

\define\nl{\newline}
\redefine\i{^{-1}}

\define\un{\underline}

\define\ot{\otimes}
\define\bbq{\bar{\QQ}_l}

\define\End{\text{\rm End}}

\define\tr{\text{\rm tr}}

\define\ph{\phi}

\define\r{\rho}
\define\s{\sigma}

\define\kk{\bold k}

\define\qq{\bold q}

\define\CC{\bold C}

\define\FF{\bold F}

\define\NN{\bold N}

\define\QQ{\bold Q}
\define\RR{\bold R}

\define\ZZ{\bold Z}

\define\cb{\Cal B}

\define\ce{\Cal E}
\define\cf{\Cal F}

\define\ch{\Cal H}

\define\fE{\frak E}

\define\fH{\frak H}

\define\fR{\frak R}

\define\sha{\sharp}

\define\bul{\bullet}

\head Introduction\endhead
\subhead 0.1\endsubhead
Let $G$ be a connected reductive group over an algebraically closed field $\kk$. Let $W$ be the
Weyl group of $G$ and let $\cb$ be the variety of Borel subgroups of $G$.
In this paper we consider two partitions of $\cb$ into pieces indexed by the various $w\in W$.
One partition (introduced in \cite{L79}) consists of the subvarieties  

(a) $Y_{s,w}=\{B\in\cb;pos(B,sBs\i)=w\}$
\nl
where $s$ is in $G_*$ (the open subset of $G$ consisting of regular semisimple elements)
and  $pos:\cb\T\cb@>>>W$ is the relative position map. The other partition (introduced in
\cite{DL76}) is defined when

(b) $\kk$ is an algebraic closure of a finite prime field;
\nl
and

(c) $F:G@>>>G$ is the Frobenius map for an $\FF_q$-split rational structure on $G$
\nl
(with $\FF_q$ being the subfield of $\kk$ with $q$ elements); it consists of the subvarieties 

(d) $X_w=\{B\in\cb;pos(B,F(B))=w\}$.
\nl
(The definition of the varieties in (a) was inspired by that of the varieties in (d).)
One of the themes of this paper is to point out a remarkable similarity between $Y_{s,w}$ in (a) and the
quotient $U^F\bsl X_w$ of $X_w$ in (b), where $U^F$ is the group of rational points of the
unipotent radical of an $F$-stable $B\in\cb$, acting by conjugation.

We show that $Y_{s,w}$ is affine when $w$ has minimal length in its conjugacy class (the analogous result
for $U^F\bsl X_w$ was known earlier). We show that the closure of $Y_{s,w}$ has the same type
of singularities as the closure of a Bruhat cell (the analogous result
for $X_w$ was known earlier). We show that if $\kk,F,q$ are as in (b),(c), the number of fixed point of $F^t$
on $Y_{s,w}$ and on $U^F\bsl X_w$ is the same (here $t\in\{1,2,\do\}$); this result is implicit
in \cite{L78}, \cite{L79}. (This number can be expressed as a trace of left multiplication by the standard basis
element $T_w$ on the Iwahori-Hecke algebra of $W$.) We show that while the cohomologies
of $X_w$ give rise to a virtual $G^F$-module (see \cite{DL76}) (so that the cohomologies
of $U^F\bsl X_w$ give rise to a virtual module of the Hecke algebra of $G^F$ with respect to $B^F$),
the cohomologies of $Y_{s,w}$ give rise to a virtual $W$-module.

But there is one key difference between $Y_{s,w}$ and $U^F\bsl X_w$: the former can be defined also
in characteristic zero, while the latter cannot. A study of $Y_{s,w}$ over the real numbers
can be found in \S3.

Another
biproduct of our study is a way to associate to any $w$ in $W$ (or more generally any Coxeter group) a subset
$\fE(w)$ of $W$, see 1.10 and \S5.

\subhead 0.2\endsubhead
For any $B\in \cb$ ket $U_B$ be the unipotent radical of $B$. For $B\in\cb,w\in W$ let
$$\cb_{B,w}=\{B'\in\cb;pos(B,B')=w\},$$
$$\cb_{B,\le w}=\{B'\in\cb;pos(B,B')\le w\}.$$
Here $\le$ is the standard partial order on $W$.
For $w$ in $W$ and $B,B'$ in $\cb$ let
$$Z_{B,B',w}=\{B''\in\cb;pos(B,B'')=z,pos(B',B'')=w\},$$
$$Z_{B,B',\le w}=\{B''\in\cb;pos(B,B'')=z,pos(B',B'')\le w\}$$
where $z=\pos(B,B')$.
For a maximal torus $T$ of $G$ we set $\cb^T=\{B\in\cb;T\sub B\}$. 
For $s\in G_*$ let $T_s$ be the unique maximal torus containing $s$.

\proclaim{Proposition 0.3} Let $s\in G_*,w\in W,z\in W$.
Let $B\in\cb^{T_s},B'\in\cb^{T_s}$ be such that $pos(B,B')=z$.

(a) There is a canonical isomorphism $Y_{s,w}\cap\cb_{B,z}@>\si>>Z_{B,B',w}$.

(b) Let $Y_{s,\le w}=\cup_{y\in W}Y_{s,y}$ (a closed subset of $\cb$).
There is a canonical isomorphism $Y_{s,\le w}\cap\cb_{B,z}@>\si>>Z_{B,B',\le w}$.
\endproclaim
Statement (a) is contained in the proof of \cite{L79, 1.2}.

\proclaim{Proposition 0.4} Assume that $\kk,F$ are as in 0.1(b),(c). Let $w\in W,z\in W$.
Let $B\in\cb^F,B'\in\cb^F$ be such that $pos(B,B')=z$.

(a) There is a canonical isomorphism $U_B^F\bsl(X_w\cap\cb_{B,z})@>\si>>Z_{B,B',w}$.

(b) Let $X_{\le w}=\cup_{y\in W}X_y$ (a closed subset of $\cb$).
There is a canonical isomorphism $U_B^F\bsl(X_{\le w}\cap\cb_{B,z})@>\si>>Z_{B,B',\le w}$.
\endproclaim
The action of $U_B^F$ is by conjugation.

\subhead 0.5\endsubhead
Let $\qq$ be an indeterminate. Let $\ch$ be the free $\ZZ[\qq]$-module with basis
$\{T_w;w\in W\}$. It is well known that there is a unique structure of associative
$\ZZ[\qq]$-algebra on $\ch$ such that $T_yT_{y'}=T_{yy'}$ if $l(yy')=l(y)+l(y')$
and $(T_y+1)(T_y-\qq)=0$ if $l(y)=1$. Here $l:W@>>>\NN$ is the length function.

For $w,w'$ in $W$ we have
$T_wT_{w'}=\sum_{w''\in W}N_{w,w',w''}T_{w''}$ where $N_{w,w',w''}\in\ZZ[\qq]$.

\proclaim{Corollary 0.6} Assume that $\kk,F,q$ are as in 0.1(b),(c). Let $w\in W,z\in W$.

(a) Assume that $s\in G_*\cap G^F$. Let $B\in\cb^{T_s}\cap B\in\cb^F$.
Let $t\in\{1,2,\do\}$. Then $\sha(Y_{s,w}\cap\cb_{B,z})^{F^t}=N_{w,z\i,z\i}(q^t)$.

(b) Let $B\in\cb^F$. Then $\sha(U_B^F\bsl(X_w\cap\cb_{B,z}))^{F^t}=N_{w,z\i,z\i}(q^t)$.
\endproclaim

\proclaim{Corollary 0.7} Assume that $\kk,F,q$ are as in 0.1(b),(c) and $t\in\{1,2,\do\}$.
Let $w\in W$.

(a) Assume that $s\in G_*\cap G^F$ and that $T_s$ is split over $F_q$. Then
$\sha(Y_{s,w}^{F^t})=\sum_{z\in W}N_{w,z\i,z\i}(q^t)=\tr(T_w:\ch@>>>\ch)_{\qq=q^t}$.

(b) Let $B\in\cb^F$. Then $\sha(U_B^F\bsl X_w)^{F^t}=\sum_{z\in W}N_{w,z\i,z\i}(q^t)
=\tr(T_w:\ch@>>>\ch)_{\qq=q^t}$.
\endproclaim
Here $T_w:\ch@>>>\ch$ is left multiplication by $T_w$ in $\ch$.
Note that (a) can be deduced from \cite{L79, 1.2}; (b) appears in \cite{L78, 3.10(a)}.

\proclaim{Proposition 0.8} Let $w\in W,i\in\ZZ$.  Assume that $\kk,F,q$ are as in 0.1(b),(c).

(a) If $s\in G_*\cap G^F$ and some/any $B\in\cb^{T_s}$ satisfies $F(B)=B$,
then any eigenvalue of $F$ on $H^i_c(Y_{s,w})$ is in $\{q^j;j\in\ZZ\}$.

(b) If $B\in\cb^F$, any eigenvalue of $F$ on $H^i_c(U_B^F\bsl X_w)$ is in $\{q^j;j\in\ZZ\}$.
\endproclaim

\subhead 0.9\endsubhead
Let $\un W$ be the set of conjugacy classes in $W$.
For $w\in W$ we denote by $\un w$ the conjugacy class of $w$.
In this subsection $\kk$ is as in 0.1(b). Assuming that $F$ is as in
0.1(c) we define a map $G_*\cap G^F@>>>\un W$, $s\m[s]$ by $pos(B,F(B))\in[s]$
for some/any $B\in\cb^{T_s}$. Let $w\in W$. As shown in \cite{DL76},
if $F$ is as in 0.1(c), the finite group group $G^F$ acts naturally on $H^i_c(X_w)$
($\bbq$-cohomology with compact support, $i\in\ZZ$), giving rise to a virtual representation
$R_w=\sum_i(-1)^iH^i_c(X_w)$ of $G^F$; moreover, $R_w$ depends only on $\un w$.
\nl
For $j\in\ZZ$ we denote by $H^i_c(X_w)_j$ the part of weight $j$ of $H^i_c(X_w)$ and we set
$R_{j,w}=\sum_i(-1)^iH^i_c(X_w)_j$ (a virtual representation of $G^F$). Note that $R_w=\sum_jR_{j,w}$.

The following result is proved in \S2.

(a) {\it For any $j\in\ZZ$ there is a unique virtual representation
$\fR_{j,w}$ of $W$ such that
for any $F,q$ as in 0.1(c) and any $s\in G_*\cap G^F$ we have
$\tr(y,\fR_{j,w})=\sum_{i\in\ZZ}(-1)^i\tr(F,H^i_c(Y_{s,w})_j)q^{-j/2}$ for
some/any $y\in[s]$. Here $H^i_c(Y_{s,w})_j$ is the part of weight $j$ of $H^i_c(Y_{s,w})$.
Moreover, $\fR_w:=\sum_j\fR_{j,w}$ depends only on $\un w$. }
\nl
Another definiton of $\fR_{j,w}$ (valid also in characteristic zero) is given in \cite{L79, 1.2}
(where it is denoted by $\r_{j,w}$). We will not attempt to compare it with the present definition.
In \cite{L79} several arguments are based on
a statement in \cite{L79, p.327, line -6} which was stated without proof. That statement can be proved
using the theory of character sheaves; this will not be done here.

\subhead 0.10\endsubhead
Let $w\in W$. Let
$$\fE(w)=\{z\in W; N_{w,z,z}\ne0\}.\tag a$$
From 0.6 we see that the following three conditions for $z\in W$ are equivalent.

(i) We have
$Y_{s,w}\cap\cb_{B,z\i}\ne\emp$ for some, or equivalently, any $s\in G_*,B\in\cb^{T_s}$.

(ii) We have $X_w\cap\cb_{B_0,z\i}\ne\emp$ for some, or equivalently, any $B\in\cb^F$ (here $\kk,F$
are as in 0.1(b),(c)).

(iii) We have $z\in\fE(w)$.
\nl
We use that $N(w,z,z)\ne0$ implies $N(w,z,z)(q)\ne0$ for any $q\in\{2,3,\do\}$, since

(b) $N(w,w',w'')$ is a polynomial in $\qq,\qq-1$ with coefficients in $\NN$.
\nl
I believe that the subsets $\fE(w)$ of $W$ deserve further study. A beginning of such a study can be found in \S5.

\subhead 0.11\endsubhead
I thank Xuhua He for useful discussions.

\head 1. Proof of Propositions 0.3, 0.4, 0.8\endhead
\subhead 1.1\endsubhead
We prove Proposition 0.3(a). There is a unique $B'{}^-\in\cb^{T_s}$ such that $B'\cap B'{}^-=T_s$.
Let $U_z=U_B\cap U_{B'{}^-}$. We identify $U_z$ with $\cb_{B,z}$
by $u\m uB'u\i$. Hence we can identify
$$\align&Y_{s,w}\cap\cb_{B,z}=\{u\in U_z;pos(uB'u\i,suB'u\i s\i)=w\}\\&=
\{u\in U_z;pos(B',u\i suB'u\i s\i u)=w\}\\&=\{u\in U_z;pos(B',s\i u\i suB'u\i s\i us)=w\}.\endalign$$
Using the isomorphism $U_z@>\si>>U_z$ given by $u\m s\i u\i su$ (recall that $s\in G_*$) we obtain an
identification      

(a) $Y_{s,w}\cap\cb_{B,z}=\{u'\in U_z; pos(B',u'B'u'{}\i)=w\}$.
\nl
Now $u'\m u'B'u'{}\i$ identifies

$\{u'\in U_z; pos(B',u'B'u'{}\i)=w\}=\{B''\in\cb_{B,z};pos(B',B'')=w\}=Z_{B,B',w}$.
\nl
Combining with (a) we obtain 0.3(a). The same proof (replacing $=w$ by $\le w$) gives 0.3(b).
Note that the isomorphisms in 0.3 are given by $uB'u\i\m s\i u\i su B'u\i s\i us$ with $u\in U_z$.

\subhead 1.2\endsubhead
We prove Proposition 0.4(a). Let $B'{}^-\in\cb^F$ such that
$B'\cap B'{}^-$ is a maximal torus contained in $B$.
Let $U_z=U_B\cap U_{B'{}^-}$. We identify $U_z$ with $\cb_{B,z}$
by $u\m uB'u\i$. Hence we can identify
$$\align& X_w\cap\cb_{B,z}=\{u\in U_z;pos(uB'u\i,F(uB'u\i))=w\}\\&=
\{u\in U_z;pos(B',u\i F(u)B'F(u\i)u)=w\}.\endalign$$
Using the isomorphism $U_s^F\bsl U_z@>\si>>U_z$ given by $u\m u\i F(u)$ (coming from Lang's theorem)
 we obtain an identification
$$U_z^F\bsl(X_w\cap\cb_{B,z})=\{u'\in U_z;pos(B',u'B'u'{}\i)=w\}.$$
This can be identified with $Z_{B,B',w}$ as in 1.1. This proves 0.4(a).
The same proof (replacing $=w$ by $\le w$) gives 0.4(b).
Note that the isomorphisms in 0.4 are given by $uB'u\i\m u\i F(u)B'F(u)\i u$ with $u\in U_z$.

\subhead 1.3\endsubhead
Let $T$ be a maximal torus of $G$.
Let $w_0$ be the longest element of $W$. For any $B\in\cb^T$
let ${}^B\cb=\{B'\in\cb;\pos(B',B)=w_0\}$, an open set in $\cb$. We show:

(a) $\cb=\cup_{B\in\cb^T}({}^B\cb)$.
\nl
Let $B'\in\cb$. Let $B_1\in\cb^T$. We have $pos(B_1,B')=z$ for some $z\in W$. We can find
$B\in\cb^T$ such that $pos(B,B_1)=w_0z\i$. Since $l(w_0z\i)+l(z)=l(w_0)$ we must have
$pos(B,B')=w_0$. Thus $B'\in {}^B\cb$. This proves (a).

\subhead 1.4\endsubhead
Let $w\in W,s\in G_*$. We show:

(a) {\it $Y_{s,w}$ is smooth of pure dimension $l(w)$.}
\nl
From 1.3(a) we have an open covering
$Y_{s,w}=\cup_{B\in\cb^{T_s}}(Y_{s,w}\cap{}^B\cb)$. It is enough to prove that for any
$B\in\cb^{T_s}$, $Y_{s,w}\cap{}^B\cb$ is smooth of pure dimension $l(w)$.
Define $B'\in\cb^{T_s}$ by $B\cap B'=T_s$. By 0.3(a) we can identify
$Y_{s,w}\cap{}^B\cb=Z_{B,B',w}=\{B''\in\cb;pos(B,B'')=w_0,pos(B',B'')=w\}$.
(The identification is by $uB'u\i\m s\i u\i su B'u\i s\i us$ with $u\in U_B$.)
This is the intersection of the smooth irreducible
subvariety $\{B''\in\cb;pos(B',B'')=w\}$ of dimension $l(w)$
of $\cb$ with the open subset $\{B''\in\cb;pos(B,B'')=w_0\}$ of $\cb$ hence it is
smooth of pure dimension $l(w)$. This proves (a).

Another (longer) proof of (a) is given in \cite{L79, 1.1}.

\subhead 1.5\endsubhead
In this subsection we assume that $\kk,F$ are as in 0.1(b),(c). Leet $w\in W$.
The following result appears in \cite{DL76,1.4}.

(a) {\it $X_w$ is smooth of pure dimension $l(w)$.}
\nl
We give an alternative proof of (a) similar to that in 1.4.
Let $T$ be a maximal torus of $G$ such that $F(T)=T$. From 1.3(a) we have an open covering
$X_w=\cup_{B\in\cb^T}(X_w\cap{}^B\cb)$. It is enough to prove that for any
$B\in\cb^T$, $X_w\cap{}^B\cb$ is smooth of pure dimension $l(w)$.

Define $B'\in\cb^T$ by $B\cap B'=T$. By 0.4(a) we can identify
$U_B^F\bsl(X_w\cap{}^B\cb)=Z_{B,B',w}$ where as in the proof in 1.4, 
$Z_{B,B',w}$ is smooth of pure dimension $l(w)$.
The conjugation action of $U_B^F$ on $X_w\cap {}^B\cb$ is free since the conjugation action of
$U_B$ on ${}^B\cb$ is free (and transitive). It follows that
$X_w\cap{}^B\cb$ is smooth of pure dimension $l(w)$. This proves (a).

\subhead 1.6\endsubhead
Let $s\in G_*,w\in W$. Let $T=T_s$. We show:

(a) $Y_{s,\le w}$ is equal to the closure of $Y_{s,w}$ in $\cb$.
\nl
Let $B_1\in Y_{s,\le w}$. We have $B_1\in Y_{s,y}$ for a unique $y\in W$.
By 1.3(a) we can find $B\in\cb^T$ such that $B_1\in{}^B\cb$.
Define $B'\in\cb^T$ by $B\cap B'=T$.
By 0.3(b) the open set $Y_{s,\le w}\cap{}^B\cb$ of $Y_{s,\le w}$ is identified
${}^B\cb\cap\cb_{B',\le w}$, an open set in $\cb_{B',\le w}$.
Under this identification $B_1$ becomes an element $B_2\in{}^B\cb\cap\cb_{B',y}$.
Since ${}^B\cb$ is open in $\cb$ we have

${}^B\cb\cap(\text{ closure of $\cb_{B',w}$ in }\cb)\sub
\text{ closure of ${}^B\cb\cap\cb_{B',w}$ in}{}^B\cb$.
\nl
In particular we have $B_2\sub\text{ closure of ${}^B\cb\cap\cb_{B',w}$ in}{}^B\cb$.
Using again the identification above we deduce that $B_1$ is in the closure of
$Y_{s,w}\cap {}^B\cb$ in ${}^B\cb$ and in particular $B_1$ is in the closure of
$Y_{s,w}$ in $\cb$. This proves (a).

Now let $B_1,B,B',B_2,y$ be as in the proof of (a).
Let $\fH^i_{B_1}$ be the stalk at $B_1$ of the $i$-th cohomology sheaf of
the intersection cohomology complex of $Y_{s,\le w}$ with coefficients in $\bbq|Y_{s,w}$
(this is defined in view of (a) and 1.4(a)). From the proof of (a) we see that

(b) $\fH^i_{B_1}={}'\fH^i_{B_2}$,
\nl
where ${}'\fH^i_{B_2}$ is the stalk at $B_2\in\cb_{B',y}$ of the $i$-th cohomology sheaf of
the intersection cohomology complex of $\cb_{B',\le w}$ with coefficients in $\bbq|\cb_{B',w}$.

Note that if $w\in W$ and $\kk,F$ are as in 0.1(b),(c), then results similar to (a),(b)
are known to hold for $X_{\le w}$. 

\subhead 1.7\endsubhead    
We prove Proposition 0.8. Let $B$ be as in (a),(b). Using the decompositions
$$Y_{s,w}=\sqc_{z\in W}(Y_{s,w}\cap\cb_{B,z}),$$
$$U_B^F\bsl X_w=\sqc_{z\in W}(U_0^F\bsl(X_w\cap\cb_{B,z})),$$
we see that it is enough to show that for any $z\in W$ and any $i\in\ZZ$,
any eigenvalue of $F$ on $H^i_c(Y_{s,w}\cap\cb_{B,z})$ 
is in $\{q^j;j\in\ZZ\}$ (in case (a)) and any eigenvalue of $F$ on
$H^i_c(U_0^F\bsl(X_w\cap\cb_{B_0,z}))$ is in $\{q^j;j\in\ZZ\}$ (in case (b)).
Using 0.3, 0.4, we see that it is enough to show that for any $z\in W$ and any
$i\in\ZZ$, any eigenvalue of $F$ on $H^i_c(Z_{B,B',w})$ is in $\{q^j;j\in\ZZ\}$ where
$B'$ is as in 0.3, 0.4. This is a special case of \cite{L78, 3.7}.

\subhead 1.8\endsubhead    
We prove Corollary 0.6. Let $B$ be as in 0.6.  Using 0.3, 0.4 we see that it is enough
to prove that for some $B'\in\cb^F$ with $pos(B,B')=z$ we have

(a) $\sha(Z_{B,B',w}^{F^t})=N_{w,z\i,z\i}(q^t)$.
\nl
This is a special case of \cite{L78, 3.7}.

\head 2. Construction of $\fR_{j,w}$\endhead
\subhead 2.1\endsubhead
For any $n\in\ZZ$ let $\ch_n=\bbq\ot_{\ZZ[\qq]}\ch$ where $\bbq$ is viewed as a
$\ZZ[\qq]$-algebra via $\qq\m n$. Then $\ch_n$ is a $\bbq$-algebra with $\bbq$-basis
$\{T_w;w\in W\}$. If $n\ne-1$, the algebra $\ch_n$ is semisimple and the
irreducible $\ch_n$-modules (up to isomorphism) are in natural bijection
$E_n\lra E$ with $\Irr(W)$, the set of irreducible $W$-modules over $\bbq$
(up to isomorphism), once $\sqrt{n}$ has been chosen.

In this section we assume that $\kk$ is as in 0.1(b). Let $w\in W,j\in\ZZ$.
Note that the uniqueness of $\fR_{j,w}$ in 0.9(a) is clear since for any $y\in W$
we can find $F,q$ as in 0.1(c) and $s\in G_*\cap G^F$ such that $y\in[s]$. We now prove
the existence of $\fR_{j,w}$ in 0.9(a).

Let $F,q$ be as in 0.1(c). Let $\cf_q$ be the vector space of functions $\cb^F@>>>\bbq$.
By assigning to the basis element $T_w$ of $\ch_q$ 
the linear map $T_w:\cf_q@>>>\cf_q$ given by $f\m f'$ 
where $f'(B)=\sum_{B'\in\cb^F;pos(B,B')=w}f(B')$, we identify $\ch_q$ with
a subalgebra of $\End(\cf_q)$. We have a canonical
decomposition $\cf_q=\op_{E\in\Irr(W)}E_q\ot[E]_q$ (as a $(\ch_q,G^F)$-module) where
$[E]_q$ is an irreducible representation of $G^F$.

Let $s\in G_*\cap G^F$. Let $|W|$ be the order of $W$.
For $t\in\NN$ we set $F_t=F^{1+t|W|},q_t=q^{1+t|W|}$. We have
$$\align&\sum_{i\in\ZZ}(-1)^i\tr(F_t,H^i_c(Y_{s,w}))\\&=\sha(B\in Y_{s,w};F_t(B)=B)=
\sha(B\in\cb^{F_t};pos(B,sBs\i)=w)\\&=
\tr(sT_w:\cf_{q_t}@>>>\cf_{q_t})=\sum_{E\in\Irr(W)}\tr(T_w,E_{q_t})\tr(s,[E]_{q_t}).
\endalign$$
By \cite{DL76, 7.8} we have $\tr(s,[E]_{q_t})=([E]_{q_t}:R_{y,q_t})$ where
$y\in[s]$, $R_{y,q_t}$ is defined as $R_y$ in 0.3(a) with $F,q$ replaced by $F_t,q_t$
and $([E]_{q_t}:R_{y,q_t})$ denotes multiplicity of an irreducible $G^{F_t}$-module
in a virtual $G^{F_t}$-module. Note that the conjugacy class $[s]$ in $W$ associated
to $s$ and to $F_t$ is independent of $t$ since $F^{|W|}(B)=B$ for any $B\in\cb$
that contains $s$.
Moreover, from \cite{L84} it is known that $([E]_{q_t}:R_{y,q_t})$ is independent of $t$.
Thus we have
$$\sum_{i\in\ZZ}(-1)^i\tr(F_t,H^i_c(Y_{s,w}))=
\sum_{E\in\Irr(W)}\tr(T_w,E_{q_t})([E]_q:R_y).$$
From 0.8 we see that $H^i_c(Y_{s,w})_{j'}=0$ if $j'$ is odd, while if $j'$ is even
$F^{|W|}$ acts on $H^i_c(Y_{s,w})_{j'}$ as $q^{j'|W|/2}$ times a unipotent transformation
so that $\tr(F_t,H^i_c(Y_{s,w})_{j'})=\tr(F,H^i_c(Y_{s,w})_{j'})q^{j't|W|/2}$. We see that
$$\sum_{i,j'\text{ in }\ZZ}(-1)^i\tr(F,H^i_c(Y_{s,w})_{j'})q^{-j'/2}
q_t^{j'/2}=\sum_{E\in\Irr(W)}\tr(T_w,E_{q_t})([E]_q:R_y).$$
It is known that for $E\in\Irr(W)$ we have
$\tr(T_w,E_{q_t})=\sum_{j'\in\ZZ}\mu_{w,E;j'}q_t^{j'/2}$
for $t\in\NN$ where $\mu_{w,E;j'}\in\ZZ$ are independent of $t$ and are zero for all
but finitely many $j'$. Thus we have
$$\sum_{i,j'\text{ in }\ZZ}(-1)^i\tr(F,H^i_c(Y_{s,w})_{j'})q^{-j'/2}
q_t^{j'/2}=\sum_{E\in\Irr(W),j'\in\ZZ}\mu_{w,E;j'}([E]_q:R_y)q_t^{j'/2}\tag a$$
for $t\in\NN$.
By comparing the coefficient of $q_t^{j/2}$ in the two sides of (a), we obtain
$$\sum_{i\in\ZZ}(-1)^i\tr(F,H^i_c(Y_{s,w})_j)q^{-j/2}
=\sum_{E\in\Irr(W)}\mu_{w,E;j}([E]_q:R_y).\tag b$$
Now for $E\in\Irr(W)$ we set $R_E=|W|\i\sum_{z\in W}tr(z,E)R_z$
(a rational linear combination of representations of $G^F$)
and we use the equality $([E]_q:R_{E'})=([E']_q:R_E)$ for $E,E'$ in $\Irr(W)$
(a known property \cite{L84} of the nonabelian Fourier transform).
We obtain
$$\align&\sum_{i\in\ZZ}(-1)^i\tr(F,H^i_c(Y_{s,w})_j)q^{-j/2}
=\sum_{E,E'\text{ in }\Irr(W)}\tr(y,E')([E]_q:R_{E'})\mu_{w,E;j}\\&
=\sum_{E,E'\text{ in }\Irr(W)}\tr(y,E')([E']_q:R_E)\mu_{w;E;j}\\&
=\sum_{E'\in\Irr(W)}\tr(y,E')([E']_q:\sum_{E\in\Irr(W)}\mu_{w,E;j}R_E)
=\sum_{E'\in\Irr(W)}\tr(y,E')([E']_q: R_{j,w})\endalign$$
where the last equality can be deduced from \cite{L84, 3.8}.
We see that $\fR_{j,w}:=\sum_{E'\in\Irr(W)}([E']_q:R_{j,w})E'$
has the properties stated in 0.9(a).
Note that $\fR_w=\sum_j\fR_{j,w}=\sum_{E'\in\Irr(W)}([E']_q:R_w)E'$
depends only on $\un w$ since $R_w$ has such a property. This completes the proof of 0.9(a).

\subhead 2.2\endsubhead
From the proof in 2.1 we see that for $E\in\Irr(W),w\in W,j\in2\NN$ we have

(a) $(E:\fR_{j,w})=([E]_q:R_{j,w})$
\nl
where the left hand side is a multiplicity as a $W$-module
and the right hand side is a multiplicity as a $G^F$-module. It follows that 

(b) $(E:\fR_w)=([E]_q:R_w)$.

\head 3. Over real numbers\endhead
\subhead 3.1\endsubhead
In this section we assume that $\kk=\CC$ and that $G$ has a given split $\RR$-structure.
Let $G(\RR)$, $\cb(\RR)$ be the set of real points of $G,\cb$.
Let $s\in G_*\cap G(\RR)$ be such that for any $B\in\cb^{T_s}$ we have $B\in\cb(\RR)$.
Let $w\in W$ and let $Y_{s,w}(\RR)=Y_{s,w}\cap\cb(\RR)$. This is a smooth manifold
of pure (real) dimension equal to $l(w)$. (See 1.4(a).)

For $w$ in $W$ and $B,B'$ in $\cb^{T_s}$ let
$Z_{B,B',w}$ be as in 0.2. Then $Z_{B,B',w}$ is defined over $\RR$ and we set
$Z_{B,B',w}(\RR)=Z_{B,B',w}\cap\cb(\RR)$. Let $z=\pos(B,B')$.
The following result can be deduced from \cite{R98, 6.1}:

(a) $\chi_c(Z_{B,B',w}(\RR))=N_{w,z\i,z\i}(-1)$
\nl
where $\chi_c$ is Euler characteristic in cohomology with compact support.
Using (a) and the homeomorphism $Y_{s,w}(\RR)\cap\cb_{B,z}\cong Z_{B,B',w}(\RR)$ deduced from
0.3(a) we see that 

$\chi_c(Y_{s,w}(\RR)\cap\cb_{B,z})=N_{w,z\i,z\i}(-1)$.
\nl
Using this and the partition $Y_{s,w}(\RR)=\cup_{z\in W}(Y_{s,w}(\RR)\cap\cb_{B,z})$,
we deduce:

\proclaim{Proposition 3.2} We have
$\chi_c(Y_{s,w}(\RR))=\sum_{z\in W}N_{w,z\i,z\i}(-1)=\tr(T_w:\ch_{-1}@>>>\ch_{-1})$.
\endproclaim
Here $T_w:\ch_{-1}@>>>\ch_{-1}$ is left multiplication by $T_w$ in $\ch_{-1}$.

\head 4. Affineness\endhead
\subhead 4.1\endsubhead
Let $s\in G_*,w\in W$. Assume that $w$ has minimal length in $\un w$. We show:

(a) {\it $Y_{s,w}$ is affine.}
\nl
The proof is a modification of the proof of the analogous result for $X_w$ given in
\cite{BR08}. As in \cite{loc.cit.} we can assume that $w$ is elliptic and good (as defined
in \cite{loc.cit.}). For such $w$, it is shown in \cite{loc.cit.} that for some $n\ge1$,
the variety
$$\align&V=\{(B_0,B_1,B_2,\do,B_n)\in\cb^{n+1};\\&pos(B_0,B_1)=pos(B_1,B_2)=\do=pos(B_{n-1},B_n)=w\}\endalign$$
is affine. Define $\ph:\cb@>>>\cb^{n+1}$ by $B\m(B,sBs\i,s^2Bs^{-2},\do,s^nBs^{-n})$.
Then $\ph$ is an isomorphism of $\cb$ onto a closed subvariety $V'$ of $\cb^{n+1}$.
Then $V\cap V'$ is closed in $V$ hence is affine. Now $\ph$ restricts to an isomorphism
of $Y_{s,w}$ onto $V\cap V'$ hence $Y_{s,w}$ is affine.

\head 5. The subset $\fE(w)$ of $W$\endhead
\subhead 5.1\endsubhead
In this section $W$ is allowed to be any Coxeter group. For any $w,w',w''$ in $W$, the
polynomials $N_{w,w',w''}\in\ZZ[\qq]$ can be defined in terms of the Iwahori-Hecke algebra $\ch$
with basis $\{T_w;w\in W\}$ attached as in 0.5 to $W$. As in 0.10(a) for any $w\in W$ we define a
subset of $W$ by
$$\fE(w)=\{z\in W;N_{w,z,z}\ne0\}.$$
For example, if $w=1$, then $\fE(w)=W$; if $w=\s$, $l(\s)=1$,
then $\fE(w)=\{z\in W;l(\s z)=l(z)-1\}$. 

(a) {\it If $W$ is finite and $w_0$ is the longest element of $W$, then for any $w$ we have
$w_0\in\ce(w)$. In particular, $\fE(w)\ne\emp$.}
\nl
We argue by induction on $l(w)$. If $w=1$ we have $T_wT_{w_0}=T_{w_0}$ and the desired result
holds. Assume now that $l(w)>0$. We have $w=\s z$ where $l(\s)=1,l(z)=l(w)-1$. From the
definition we have $N_{w,w_0,w_0}=N_{z,w_0,w_0}(\qq-1)+N_{z,w_0,\s w_0}$. Since any $N_{a,b,c}$ is a
satisfies 0.10(b) and $N_{z,w_0,w_0}\ne0$ by the induction
hypothesis, it follows that $N_{w,w_0,w_0}\ne0$. This proves (a).

(b) {\it If $W$ is a finite Weyl group and $w$ is a Coxeter element of minimal length in $W$
then $\fE(w)=\{w_0\}$.}
\nl
Assume that $\kk,F,q$ are as in 0.1(b),(c). Let $B,B'$ be in $\cb^F$.
Using \cite{L76, 2.5} we have $X_w\cap\cb_{B,z}=\emp$ for $z\ne w_0$.
Using this and 0.4(a), we see that $Z_{B,B',w}=\emp$ if $pos(B,B')\ne w_0$.
Using now 0.6, we see that
$N_{w,z\i,z\i}(q)=0$ if $z\ne w_0$. Since $q$ can take infinitely many values we see
that $N_{w,z\i,z\i}=0$ if $z\ne w_0$, so that $\fE(w)\sub\{w_0\}$. By (a) we have
$\fE\ne\emp$ and (b) follows.

(c) {\it Assume that $n\in\{2,3,...\}$ and that $W$ has generators $\s_1,\s_2$ and relations
$\s_1^2=\s_2^2=1$ and $(\s_1\s_2)^{2n}=1$ (a dihedral group of order $4n$).
Let $w=(s_1s_2)^k$ where $k\in\{1,2,\do,n\}$. We have $\fE(w)=\{z\in W;l(z)\ge2n-k+1\}$.}
\nl
The proof is by computation.

(d) {\it Assume
that $W$ has generators $\s_1,\s_2$ and relations $\s_1^2=\s_2^2=1$ (an infinite
dihedral group). Let $w=(s_1s_2)^k$ where $k\in\{1,2,\do\}$. We have $\fE(w)=\emp$.
We have $\fE(s_1s_2s_1)=\{s_1s_2,s_1s_2s_1,s_1s_2s_1s_2,s_1s_2s_1s_2s_1,\do\}$.}
\nl
The proof is by computation.

\subhead 5.2\endsubhead
Let $W^\bul=\{w\in W;\fE(w)\ne\emp\}$.
If $W$ is finite then by 5.1(a) we have $w_0\in\fE(w)$ for any $w\in W$; thus in this case $W^\bul=W$.
If $W$ is infinite then it may happen that $W^\bul\ne W$. For example in the setup of 5.1(d) we have
$s_1s_2\n W^\bul$.

(a) {\it For any $w\in W^\bul$ and any $z\in\fE(w)$ we have $\deg(N_{w,z,z})\le l(w)$.}
\nl
The proof is immediate; see for example \cite{L20, 2(b)}.

For any $w\in W^\bul$ we set

(b) $d(w)=\max_{z\in \fE(w)}\deg N_{w,z,z}$.
\nl
From (a) we have

(c) $d(w)\le l(w)$ for all $w\in W^\bul$.
\nl
If $W$ is finite we have $d(w)=l(w)$ for all $w\in W=W^\bul$; indeed, by the proof of 5.1(a), we have 
$\deg N_{w,w_0,w_0}=l(w)$.  If $W$ is infinite  then it may happen that $d(w)<l(w)$ for some $w\in W^\bul$.
For example in the setup of 5.1(d) we have $d(s_1s_2s_1)=2<3=l(s_1s_2s_1)$.

For $w\in W^\bul$ we define $\fE'(w)=\{z\in\fE(w);\deg N_{w,z,z}=d(w)\}$ (a nonempty set).
If $W$ is finite and $w$ is such that any simple reflection of $W$ appears in any reduced decomposition of $w$
then

(d) $\fE'(w)=\{w_0\}$.
\nl
This can be deduced from \cite{L78, p.29, lines 2-4} (which was stated without proof) or it can be deduced
from results in \cite{L20, no.2}.

\subhead 5.3\endsubhead
We return to the setup in 0.1 and we assume that $G$ is adjoint, that
$s\in G_*$ and that $w\in W$ is elliptic, of minimal length in $\un w$. 
Then $T=T_s$ acts on $Y_{s,w}$ by conjugation. The following result can be deduced from \cite{L11, 5.2}.

(a) {\it Any isotropy group of the $T$-action on $Y_{s,w}$ is finite.}
\nl
In the case where $w$ is a Coxeter element of minimal length, we have the
following stronger result.

(b) {\it $Y_{s,w}$ is a principal homogeneous space for $T$.}
\nl
Let $B\in\cb^T$. Let $B'\in\cb^T$ be such that $B\cap B'=T$.
From 5.1(b) we see (using 0.6) that $Y_{s,w}\cap\cb_{B,z}=\emp$ if $z\ne w_0$. Hence
$Y_{s,w}=Y_{s,w}\cap\cb_{B,w_0}$. From 1.1(a) we can identify $Y_{s,w}=Y_{s,w}\cap\cb_{B,w_0}$
(with its $T$-action by conjugation) and $\{u'\in U_B; pos(B',u'B'u'{}\i)=w\}$
(with its $T$-action by conjugation); by \cite{L76, 2.2} this last $T$-space is principal homogeneous.
This proves (b).

\enddocument